\documentclass[12pt, reqno]{amsart}

\usepackage[margin=1.1in]{geometry}

\usepackage{amsfonts}
\usepackage{amssymb}
\usepackage{amsmath}
\usepackage{amsthm}
\usepackage{latexsym}
\usepackage{graphicx}
\usepackage{xypic}
\usepackage{hyperref}
\usepackage{float}
\usepackage{soul}

\hypersetup{colorlinks,citecolor=red,pdfstartview=FitH, pdfpagemode=None}

\numberwithin{equation}{section}

\theoremstyle{definition}
\newtheorem{theorem}{Theorem}[section]
\newtheorem{lemma}[theorem]{Lemma}
\newtheorem{proposition}[theorem]{Proposition}

\newtheorem{remark}[theorem]{Remark}

\renewcommand{\ge}{\geqslant}
\renewcommand{\le}{\leqslant}

\def\GL{{\rm GL}}

\font\germ=eufm10

\def\a{\mbox{\germ a}}

\def\g{\mbox{\germ g}}

\def\k{\mbox{\germ k}}

\def\p{\mbox{\germ p}}

\def\R{\mathbb R}
\def\C{\mathbb C}

\def\M{\mathcal M}
\def\H{\mathcal H}
\def\P{\mathcal P}

\def\U{\mathcal U}

\def\tr{{\rm tr\,}}

\def\ad{{\operatorname{ad}\,}}
\def\Ad{{\operatorname{Ad}\,}}

\def\conv{{\mbox{conv}\,}}

\begin{document}

\title[On two geometric means and sum of adjoint orbits]{On two geometric means and sum of adjoint orbits}

\author{Luyining Gan}
\address{Department of Mathematics and Statistics\\ University of Nevada, Reno\\ Reno \\ NV 89557-0084\\ USA}
\email{lgan@unr.edu}

\author{Xuhua Liu}
\address{Department of Mathematics\\ North Greenville University \\ Tigerville \\ SC 29688 \\ USA}
\email{Roy.Liu@ngu.edu}

\author{Tin-Yau Tam}
\address{Department of Mathematics and Statistics\\ University of Nevada, Reno\\ Reno \\ NV 89557-0084\\ USA}
\email{ttam@unr.edu}

\date{\today}



\dedicatory{In honor of the retirement of Professor Ming Liao from Auburn University in 2022.} 

\begin{abstract}
{In this paper,
we study the metric geometric mean introduced by Pusz and Woronowicz and the spectral geometric mean introduced by Fiedler and Pt\'ak, originally  for  positive definite matrices. The relation between $t$-metric geometric mean and $t$-spectral geometric mean is established via log majorization. 
The result  is then extended in the context of symmetric space associated with a noncompact semisimple Lie group. For any Hermitian matrices $X$ and $Y$, So's matrix exponential formula asserts that there are unitary matrices $U$ and $V$ such that 
$$e^{X/2}e^Ye^{X/2} = e^{UXU^*+VYV^*}.$$ 
In other words, the Hermitian matrix $\log (e^{X/2}e^Ye^{X/2})$ lies in the sum of the unitary orbits of $X$ and $Y$. 
So's result is also extended to a formula for adjoint orbits associated with a noncompact semisimple Lie group.}
\end{abstract}

\maketitle

%%%%%%%%%%%%%%%%%%%%%%%%%%%%%%%%%%%%%%%%%%%%%%%%%%%%%%%
\section{Introduction}

Let $\M_n$ be the linear space of all $n \times n$ complex matrices, $\H_n$  the real subspace of $\M_n$ consisting of Hermitian matrices,  $\P_n$  the subset of positive definite matrices in $\M_n$, and $\U_n$  the group of unitary matrices in $\M_n$. For any $X \in \M_n$, let $e^X$ denote the exponential of $X$.

Given the eigenvalues of $X, Y \in \H_n$, a famous problem of Weyl \cite{We12} was to give a complete description of the eigenvalues of the matrices in the sum of the unitary orbits of $X$ and $Y$
\[
O(X) + O(Y) = :\{ UXU^{*} + VYV^{*}: U, V \in \U_n \}.
\]
This long-standing problem was completely solved. See \cite{DST98, Fu00, Kl98, KT01} and the references therein for historical development. 
The generalization of Weyl's problem to compact Lie groups was given in \cite [Theorem 9.3]{OS00}.
Its generalization to noncompact semisimple Lie groups was considered as an equivalent optimization problem in \cite{LT13}, where Liu and Tam derived the gradient flow of the smooth function associated with the optimization problem.

Let $X, Y \in \H_n$. So and Thompson \cite{ST91} conjectured and So \cite{So04}, based on a result of Klyachko \cite{Kl00},  proved that $e^{X/2}e^Ye^{X/2} = e^Z$ for some $Z \in O(X) + O(Y)$. In other words, 
\begin{equation} \label{E:So04}
\log (e^{X/2}e^Ye^{X/2}) \in O(X) + O(Y).
\end{equation}
Motivated by \eqref{E:So04}, Kim and Lim \cite{KL07} showed that 
\begin{align}
\log (e^{2X} \sharp \, e^{2Y}) &\in O(X) + O(Y), \label{E:KL07-gm} \\
\log (e^{2X} \natural \, e^{2Y}) &\in O(X) + O(Y), \label{E:KL07-sm}
\end{align}
where $\sharp$ and $\natural$ denote the metric geometric mean (geometric mean for short) and spectral geometric mean (spectral mean for short), respectively, defined for $A, B \in \P_n$ as
\begin{align*}
A \sharp B &= A^{1/2}(A^{-1/2}BA^{-1/2})^{1/2}A^{1/2}, \\
A \natural B &= (A^{-1}\sharp B)^{1/2}A(A^{-1}\sharp B)^{1/2}.
\end{align*}

The geometric mean was introduced by Pusz and Woronowicz \cite {PW75}, while the spectral mean by Fiedler and Pt\'ak \cite {FP97}. Fiedler and Pt\'ak named them as metric geometric mean and spectral geometric mean, respectively, because $(A \natural B)^2$ is similar to $AB$ and thus the eigenvalues of $A \natural B$ are the positive square roots of the corresponding eigenvalues of $AB$. 
It is known that $A \sharp B = B \sharp A$ and $A \natural B = B \natural A$.
When $A$ and $B$ commute, we have $A \sharp B  = A \natural B = A^{1/2}B^{1/2}$.

For $A, B \in \P_n$ and $t \in [0, 1]$, the $t$-metric geometric mean ($t$-geometric mean for short) and $t$-spectral geometric mean ($t$-spectral mean for short) are naturally defined by
\begin{align*}
A \sharp_t B &= A^{1/2}(A^{-1/2}BA^{-1/2})^{t}A^{1/2}, \\
A \natural_t B &= (A^{-1}\sharp B)^{t}A(A^{-1}\sharp B)^{t}.
\end{align*}
Both are paths joining $A$ and $B$ in $\P_n$. 
Recently $t$-geometric mean has been gaining intensive interest, partially because of its connection with Riemannian geometry: $\P_n$ can be equipped
with a suitable Riemannian metric so that the curve $\gamma(t) = A \sharp_{t}B$ with $0 \le t\le 1$ is the unique geodesic joining $A$ and $B$ in $\P_n$. 
The $t$-geometric mean was generalized to  symmetric spaces of noncompact type \cite{LLT14}.
{The $t$-spectral mean was first introduced by Lee and Lim~\cite{LL07} in 2007. In the same year, it was also studied by Ahn, Kim and Lim \cite [p.191]{AKL07} (also see \cite [p.446]{KL07}). Its further algebraic and geometric meaning has been recently studied by Kim~\cite{Ki21}.}
When $t = 1/2$, they are abbreviated as $A \sharp_{1/2} B = A \sharp B$ and $A \natural_{1/2} B = A \natural B$.

{As the spectral  mean is defined by the geometric mean, there should be close relationship between  them. For example, Kim and Lim  \cite{KL07} gave} the following trace inequality 
\begin{equation} \label{E:KL07}
\tr (e^{rX} \sharp \, e^{rY})^{2/r} \le \tr e^{X+Y} \le \tr (e^{rX} \natural \, e^{rY})^{2/r}, \qquad \mbox{for all } \, r > 0.
\end{equation}
The first inequality in \eqref{E:KL07} was obtained by Hiai and Petz \cite{HP93} as the complement to the celebrated Golden-Thompson inequality \cite{Go65, Sy65, Th65}:
\[
\tr e^{X+Y} \le \tr e^Xe^Y.
\]
The second inequality in  \eqref{E:KL07} is actually a refinement of Golden-Thompson inequality when $0 < r \le 1$. 
See Theorem \ref{T:GT-ALT-Spectral} for a generalization of \eqref{E:KL07}.

In this paper we extend \eqref{E:So04}--\eqref{E:KL07-sm} to noncompact semisimple Lie groups. We also  study $t$-spectral mean. In particular, we show that $t$-geometric mean is log-majorized by $t$-spectral mean.
 Its generalization in the context of semisimple Lie groups is then given.
Matrix results are in Section 2, and generalizations of sum of adjoint orbits to semisimple Lie groups are in Section 3, and extensions of two geometric means for symmetric spaces are in Section 4.

%%%%%%%%%%%%%%%%%%%%%%%%%%%%%%%%%%%%%%%%%%%%%%%%%%%%%%%
\section{Matrix Inequalities}

Let $x = (x_1, x_2, \dots, x_n)$ and $y = (y_1, y_2, \dots, y_n)$ be in $\R^n$. Let $x^{\downarrow} = (x_{[1]}, x_{[2]}, \dots, x_{[n]})$ denote the rearrangement of the components of $x$ such that $x_{[1]} \ge  x_{[2]} \ge \cdots \ge x_{[n]}$. 
We say that $x$ is {\it majorized} by $y$, denoted by $x \prec y$,  if 
\[
\sum_{i=1}^k x_{[i]} \le \sum_{i=1}^k y_{[i]}, \quad k = 1, 2, \dots, n-1 \quad \text{and} \quad \sum_{i=1}^n x_{[i]} = \sum_{i=1}^n y_{[i]}.
\]
Among many equivalent conditions for majorization, the following one is suitable for generalization to Lie groups \cite{Ho54}:
\[
x \prec y \quad \Leftrightarrow \quad \conv S_n \cdot x \subset \conv S_n \cdot y,
\]
where $\conv S_n \cdot x$ denotes the convex hull of the orbit of $x$ under the action of the symmetric group $S_n$.
When $x$ and $y$ are nonnegative, we say that $x$ is {\em log-majorized} by $y$, denoted by $x \prec_{\log} y$ if 
\[
\prod_{i=1}^k x_{[i]} \le \prod_{i=1}^k y_{[i]}, \quad k = 1, 2, \dots, n-1 \quad \text{and} \quad  \prod_{i=1}^n x_{[i]} = \prod_{i=1}^n y_{[i]}.
\]
In other words,  when $x$ and $y$ are positive, $x \prec_{\log} y$ if and only if $\log x \prec \log y$, {where $\log x  = (\log x_1, \log x_2, \dots, \log x_n)$}.

For any $X \in  \M_n$, let
\[
\lambda(X) = (\lambda_1(X), \dots, \lambda_n(X)) 
\]
denote the vector of eigenvalues of $X$ whose absolute values are in non-increasing order.

According to \cite[Proposition 2.3]{KL07}, if $A, B \in \P_n$, then there exists $U \in \U_n$ such that
\begin{equation}  \label{E:KL07-sm1}
A \natural B = U(A^{1/2}BA^{1/2})^{1/2}U^*
\end{equation}
and $U(X^*X)U^* = XX^*$, where $X = (A^{-1}\sharp B)^{1/2}A^{1/2}$. In other words,  $X = (XX^*)^{1/2}U = U(X^*X)^{1/2} $, where $X$ is expressed in left and right polar decompositions. We note that such choice of $U\in \U_n$ is smooth as a  function of $A, B \in \P_n$. 
It follows from~\eqref{E:KL07-sm1} that
\[
\lambda (A \natural B) = \lambda((A^{1/2}BA^{1/2})^{1/2}) = \lambda((AB)^{1/2}).
\]
Moreover, $(A \natural B)^2$ is positively similar to $AB$ \cite[Theorem 3.2(8)]{FP97}, i.e., there exists $P \in \P_n$ such that $(A \natural B)^2 = P(AB)P^{-1}$.

The following result is a generalization of \eqref{E:KL07}. Though the statements about geometric mean have appeared in literature, we put them together for completeness and symmetry. 

\begin{theorem}  \label{T:GT-ALT-Spectral}
Let $X, Y \in \H_n$. For all $r > 0$, let 
\[
\phi (r) = (e^{rX} \sharp \, e^{rY})^{2/r} \quad \text{and} \quad \psi (r) = (e^{rX} \natural \, e^{rY})^{2/r}.
\]
Then the following statements are valid.
\begin{enumerate}
\item (Hiai and Petz 1993) $\phi (r)$ is monotonically decreasing on $r$ with respect to log-majorization:
\[
(e^{sX} \sharp \, e^{sY})^{2/s} \prec_{\log} (e^{rX} \sharp \, e^{rY})^{2/r}, \qquad \mbox{for all } \, 0 < r < s.
\]
\item $\psi (r)$ is monotonically increasing on $r$ with respect to log-majorization.
\item (Ahn, Kim and Lim 2007) $\displaystyle \lim_{r \to 0} \phi(r) = e^{X+Y} =  \lim_{r \to 0} \psi(r)$.
\item $(e^{rX} \sharp \, e^{rY})^{2/r} \prec_{\log} e^{X+Y} \prec_{\log} (e^{rX} \natural \, e^{rY})^{2/r}$ for all $r > 0$.
\end{enumerate}
In particular, $\tr (e^{rX} \sharp \, e^{rY})^{2/r}$ is monotonically decreasing on $r$ and $\tr (e^{rX} \natural \, e^{rY})^{2/r}$ is monotonically increasing on $r$. 
\end{theorem}

\begin{proof}
(1) This follows from a special case of \cite[Theorem 6.1]{TL18}.

(2) By \eqref{E:KL07-sm1}, there exists $U \in \U_n$ such that 
\[
\psi (r) = U \left( e^{rX/2}e^{rY}e^{rX/2}\right)^{1/r} U^*.
\]
Then (2) follows from \cite[Theorem 3.12]{TL18}.

(3) This is implied by the Lie-Trotter formula, as also stated in \cite[p.191]{AKL07} and \cite[p.444]{KL07}.

(4) Combining (1)--(3) yields (4).
\end{proof}

Geometric mean and $t$-geometric mean have been extensively studied, so their basic properties are well known. The following is a collection of some basic properties of spectral mean and $t$-spectral mean. \begin{proposition} \label{T:t-natural}
Let $A, B \in \P_n$ and $t \in [0, 1]$. Then
\begin{enumerate}
\item (\cite[Theorem 3.2 (4)]{FP97}) $(A \natural B)^{-1} = A^{-1} \natural B^{-1}$.
\item (\cite[p.2163]{Li12}) $(A \natural_t B)^{-1} = A^{-1} \natural_t B^{-1}$ and $A \natural_t B = B \natural_{1-t} A$.
\item (\cite[Theorem 3.1 (1)]{FP97}) $A^{-1} \sharp (A \natural B) = (A \natural B)^{-1} \sharp B$.
\item (\cite[p.2163]{Li12}) $A^{-1} \sharp (A \natural_t B) = (B \natural_t A)^{-1} \sharp B = (A^{-1} \sharp B)^t$.
\item (\cite[Theorem 3.1 (3)]{FP97}) If $C = A^{-1} \sharp (A \natural B)$, then $A \natural B = CAC = C^{-1}BC^{-1}$.
\item If $C_t = A^{-1} \sharp (A \natural_t B)$, then $A \natural_t B = C_tAC_t$ and $B \natural_t A = C_t^{-1}BC_t^{-1}$.
\item (\cite[p.2163]{Li12}) $(A \natural_r B) \natural_t (A \natural_s B) = A \natural_{(1-t)r+ts} B$ for all $r, s \in [0,1]$.
\item (\cite[Theorem 5.5 (5)]{FP97}) $A \sharp B$ is positively similar to $(A \natural B)^{1/2}U(A \natural B)^{1/2}$ for some $U \in \U_n$.
\end{enumerate}
\end{proposition}

\begin{proof}
We only need to show (6). Note that $C_t = (A^{-1} \sharp B)^t$ according to (4), so $A \natural_t B = C_tAC_t$ be definition. 
Moreover, by (2) we have 
\[
C_t^{-1}BC_t^{-1} = (A \sharp B^{-1})^tB(A \sharp B^{-1})^t =  (B^{-1} \sharp A)^t B (B^{-1} \sharp A)^t = B \natural_t A.
\qedhere\]
\end{proof}

The following proposition was mentioned by Ahn, Kim and Lim \cite[p.192]{AKL07} without proof.

\begin{proposition} \label{T:sharp-natural}
For all $A, B \in \P_n$, we have
\[
\lambda(A \sharp B) \prec_{\log} \lambda(A \natural B).
\]
\end{proposition}
\begin{proof}
By Proposition~\ref{T:t-natural} (8), there exists $U \in \U_n$ such that $A \sharp B$ is positively similar to $(A \natural B)^{1/2} U (A \natural B)^{1/2}$, which is similar to $(A \natural B)U$. By Weyl's inequality,
\[
\lambda(A \sharp B) = \lambda((A \natural B)U) = |\lambda((A \natural B)U)| \prec_{\log} s((A \natural B)U) = \lambda(A \natural B),
\]
where $s(M)$ denotes the vector of singular values of $M$ in non-increasing order.
\end{proof}

We are going to extend Proposition \ref{T:sharp-natural} from $1/2$ to $t\in [0,1]$. 
Given $A, B\in \H_n$, denote by $A\le B$ the L\"owner order, that is, $B-A$ is positive semidefinite. 
The following result is known as the joint monotonicity theorem according to Ando and Hiai \cite [p.118]{AH94}; also see \cite{Li12} in which Lim named it as L{\"{o}}wner-Heinz inequality and the original paper of L{\"{o}}wner \cite{Lo34}.
\begin{lemma}[L{\"{o}}wner-Heinz inequality]\label{lem:sharp}
For $A\le C, B\le D$ and $t\in [0,1]$,
\[A\sharp_t B \le C\sharp_t D.\]
\end{lemma}

\begin{remark} For $A, B\ge 0$, there is no such joint monotonicity result for $\prec_{\log}$. Indeed, 
 $\lambda(A\sharp B_1)\prec_{\log} \lambda(A\sharp B_2)$ is not  true in general when $\lambda(B_1) \prec_{\log} \lambda(B_2)$.
Here is a counterexample:
\[
A=\begin{bmatrix}
	16&0\\
	0&1
\end{bmatrix}, \quad
B_1 = \begin{bmatrix}
	2&0\\
	0&4
\end{bmatrix}, \quad
B_2 = \begin{bmatrix}
	1&0\\
	0&8
\end{bmatrix},
\]
satisfy $\lambda(B_1)\prec_{\log} \lambda(B_2)$
and 
\[
A\sharp B_1 = A^{1/2}B_1^{1/2} = \begin{bmatrix}
	4\sqrt{2}&0\\
	0&2
\end{bmatrix}, \quad
A\sharp B_2 = A^{1/2}B_2^{1/2} = \begin{bmatrix}
	4&0\\
	0&2\sqrt{2}
\end{bmatrix},
\]
contradicting with $\lambda(A\sharp B_1) \prec_{\log} \lambda(A\sharp B_2)$.
\end{remark}

\begin{remark}
	For $A, B\ge 0$, $A \sharp_t B \le A \natural_t B$ is not true in general for $t \in (0, 1)$. Here is a counterexample. Let
	\[
A=\begin{bmatrix}
	6&-3\\
	-3&4
\end{bmatrix} \ge 0, \quad
B = \begin{bmatrix}
	4&-2\\
	-2&5
\end{bmatrix} \ge 0.
\]
Then
\[
A\sharp B = \begin{bmatrix}
	4.8990 & -2.4495\\
	-2.4495 & 4.3870
\end{bmatrix}, \quad
A\natural B =\begin{bmatrix}
	4.8992 & -2.4896\\
	-2.4896 & 4.4273
\end{bmatrix},
\]
and the eigenvalues of $A\natural B - A\sharp B$ are $0.0651$ and $-0.0246$.
\end{remark}

Though $A \sharp_t B \le A \natural_t B$ is not true in general for $t \in (0, 1)$, we have the following result using log majorization.

\begin{theorem} \label{T:log-majorization}
For all $A, B \in \P_n$ and $t \in [0, 1]$, we have
\begin{equation}\label{eqn:main}
	\lambda(A \sharp_t B) \prec_{\log} \lambda(A \natural_t B).
\end{equation}
\end{theorem}

\begin{proof}
First, we know $\det (A\sharp_t B) = \det (A\natural_t B)$
because   
$\det (A\sharp_t B) = (\det A)^{1-t} (\det B)^t$
and
\begin{eqnarray*}
\det (A\natural_t B) &=& \det((A^{-1}\sharp B)^t A (A^{-1}\sharp B)^t)\\
& =& \det (A^{-1}\sharp B)^{2t} \det A \\
&=& [{(\det A)^{-1/2}} (\det B)^{1/2}]^{2t} \det A \\
&=& (\det A)^{1-t} (\det B)^t.
\end{eqnarray*}
Recall \cite[p.776-777]{MOA11}  that 
\[
\prod_{i=1}^{k}\lambda_i(A)=\lambda_1(C_k(A)),\quad  k=1, \dots, n,
\]
 where $C_k(A)$ denotes the $k$th compound of $A \ge 0$.
Thus, we need to show 
\[
\lambda_1(C_k(A\sharp_t B)) \le \lambda_1(C_k(A\natural_t B)), \quad k=1, \dots, n-1.
\]
Note that \cite {AH94, DAT16}
\[
	C_k(A\sharp_t B) = C_k(A)\sharp_t C_k(B), \quad k=1, \dots, n
\]
and 
\begin{eqnarray*}
	C_k(A\natural_t B) &=& C_k((A^{-1}\sharp B)^t A (A^{-1}\sharp B)^t)\\
	&=& C_k((A^{-1}\sharp B)^t)C_k(A) C_k((A^{-1}\sharp B)^t)\\
	&=& C_k(A^{-1}\sharp B)^t C_k(A) C_k(A^{-1}\sharp B)^t\\
	&=& [C_k(A)^{-1}\sharp C_k(B)]^t C_k(A) [C_k(A)^{-1}\sharp C_k(B)]^t\\
	&=& C_k(A)\natural_t C_k(B),\quad k=1, \dots, n,
\end{eqnarray*}
{where the second and third equalities are from the well-known equation $C_k(AB) = C_k(A)C_k(B)$ and the fourth equality is from $C_k(A^{-1}) = C_k(A)^{-1}$.}
So it suffices to show 
\begin{equation}\label{eqn:lambda2}
	\lambda_1(A\sharp_t B) \le \lambda_1(A\natural_t B).
\end{equation}
Note that 
\[
(\alpha A)\sharp _t (\beta B) = \alpha^{1-t}\beta^t  (A\sharp _t  B),\quad (\alpha A)\natural_t (\beta B) = \alpha^{1-t}\beta^t (A\natural_t  B), \qquad \alpha , \beta >0,
\]
that is, $A\sharp _t  B$ and $A\natural_t  B$ have the same order of homogeneity for $A, B$.
Thus we may prove that $A\natural_t B \le I$ implies $A\sharp_t B\le I$.
Next, let $C = A^{-1}\sharp B$. If $A\natural_t B \le I$, that is, $C^t A C^t \le I$, then 
\begin{equation}\label{C}A\le C^{-2t}.
\end{equation} Since 
$C = A^{-1}\sharp B,$
we have
\[
C = (A^{-1})^{1/2}  [(A^{-1})^{-1/2} B(A^{-1})^{-1/2}]^{1/2}   (A^{-1})^{1/2}.
\]
Thus by \eqref{C} and the fact that L\"owner order is invariant under congruence, we have
\begin{equation}\label{B}
B = A^{-1/2}(A^{1/2}CA^{1/2})^2A^{-1/2} = CAC \le C^{-2t+2}.
\end{equation}
Thus by Lemma~\ref{lem:sharp}, \eqref{C} and \eqref{B}, we have
\[
	A\sharp_t B \le  C^{-2t} \sharp_t C^{-2t+2}
	= (C^{-2t})^{1-t} (C^{2-2t})^{t}= I,
\]
since clearly $C^{-2t}$ and $C^{-2t+2}$ commute.
\end{proof}

%%%%%%%%%%%%%%%%%%%%%%%%%%%%%%%%%%%%%%%%%%%%%%%%%%%%%%%
\section{Sum of adjoint orbits}

We first recall some algebraic structures of semisimple Lie groups (see \cite{He78, Kn02}). 
In this section, unless otherwise specified, let $G$ be a noncompact connected semisimple Lie group with Lie algebra $\g$, let $\Theta$: $G\to G$ be a Cartan involution of $G$, and let $K$ be the fixed point set of $\Theta$, which is an analytic subgroup of $G$. Let $\theta = d\Theta$ be the differential map of $\Theta$. Then $\theta: \g \to \g$ is a Cartan involution and $\g = \k \oplus \p$ is a Cartan decomposition, where $\k$ is the eigenspace of $\theta$ corresponding to the eigenvalue $1$ (and also the Lie algebra of $K$) and $\p$ is the eigenspace of $\theta$ corresponding to the eigenvalue $-1$ (and also an $\Ad K$-invariant subspace of $\g$ complementary to $\k$). The  Killing form $B$ on $\g$ is negative definite on $\k$ and positive definite  on $\p$, and the bilinear form $B_\theta$ defined by 
\[
B_{\theta} (X, Y) = -B(X, \theta Y), \quad X, Y \in\g
\]
is an inner product on $\g$. 

For each $X \in \g$, let $e^X = \exp X$ be the exponential of $X$. Let $P = \{e^X: \, X\in\p\}$. The map $\p \times K \to G$, defined by $(X, k) \mapsto e^Xk$, is a diffeomorphism. 
So each $g \in G$ can be uniquely written as 
\begin{equation} \label{E:Cartan}
g  = pk = p(g)k(g)
\end{equation}
with $p = p(g) \in P$ and $k = k(g) \in K$. The decomposition $G = PK$ is {called} the left Cartan decomposition of $G$. Similarly, the right  Cartan decomposition is $G = KP$. 

Let $*: G \to G$ be the diffeomorphism defined by $*(g) = \Theta (g^{-1})$. We also write $g^* = *(g)$ for connivence.
Note that $*$ is not an automorphism on $G$, since $(fg)^* = g^*f^*$ for all $f, g \in G$. 
Because $K$ is the {fixed point} set of $\Theta$ and $\exp_{\g}: \p \to P$ is bijective,  we see that $p^*= p$ for all $p\in P$ and $k^*=k^{-1}$ for all $k\in K$. 
By the Cartan decomposition \eqref{E:Cartan}, we have for all $g \in G$
\begin{equation} \label{E:p}
p(g) = (gg^*)^{1/2}.
\end{equation}

Let $\a$ be a maximal abelian subspace of $\p$ and let $A$ be the analytic subgroup generated by $\a$.
Let $\a_+$ be a fixed closed Weyl chamber in $\a$ and let $A_+ = \exp \a_+$. Every element in $\p$ is $K$-conjugate to a unique element in $\a_+$. 
In other words, if $X \in \p$, there exist a unique $Z \in \a_+$ and some $k \in K$ such that 
\[
X = \Ad k (Z).
\]
We thus denote $\a_+(X) = Z$. It follows that 
\begin{equation} \label{E:G-spectral}
\exp X = \exp(\Ad k (Z)) = k\exp(Z)k^{-1} \in KA_+K.
\end{equation}
Applying \eqref{E:G-spectral} to the $P$-component of any $g \in G$ with Cartan decompositions $g = pk$, we have the following Lie group decomposition $G = KA_+K$.
In other words, each $g \in G$ can be written as
\begin{equation} \label{E:G-SVD2}
g = uav,
\end{equation}
where $u, v \in K$ and $a \in A_+$ is uniquely determined, and we denote $a_+(g) = a$.

The following result is a generalization of \eqref{E:So04} to semisimple Lie groups.

\begin{theorem} \label{T:orbit-sum-exponential}
For all $X, Y \in \p$, there exist $u, v \in K$ such that 
\begin{equation}
e^{X/2}e^Ye^{X/2} = e^{\Ad (u)X + \Ad (v)Y}.
\end{equation}
\end{theorem}

\begin{proof} 
Let $X, Y \in \p$. 
Denote $Z = \log a_+(e^{Y/2}e^{X/2}) \in \a_+$. 
Then \eqref{E:p} implies that
\[
e^{2Z} = [a_+(e^{Y/2}e^{X/2})]^2 = a_+((e^{Y/2}e^{X/2})^*(e^{Y/2}e^{X/2})) = a_+(e^{X/2}e^Ye^{X/2}).
\]
Thus there exists $w \in K$ such that 
\[
e^{X/2}e^Ye^{X/2} = we^{2Z}w^{-1} = e^{\Ad w(2Z)}.
\]
By \cite [Section 4]{AMW01} (also see \cite {EL05}), for $e^{Y/2}$ and $e^{X/2}$, there exist $Y'/2 \in \Ad K(Y/2)$ and $X'/2 \in \Ad K(X/2)$ such that 
\[
Z = Y'/2 + X'/2.
\]
Let $u', v' \in K$ be such that $X'/2 = \Ad u'(X/2)$ and $Y'/2 = \Ad v'(Y/2)$.
Then we have
\begin{align*}
e^{X/2}e^Ye^{X/2} &= e^{\Ad w (2Z)} \\
                          &= e^{\Ad w (Y' + X')} \\
                          &= e^{\Ad w (\Ad u'(X) + \Ad v'(Y))} \\
                          &= e^{\Ad (u)X + \Ad (v)Y}
\end{align*}
with $u = wu' \in K$ and $v = wv' \in K$.
\end{proof}

Let $G = PK$ be the left Cartan decomposition of $G$. The map $p \mapsto p^{1/2}K$ identify $P$ with $G/K$ as a symmetric space of noncompact type.
The $t$-geometric mean of $p, q \in P$ was defined in \cite{LLT14} as
\[
p \sharp_{t}q = p^{1/2}\left(p^{-1/2}qp^{-1/2}\right)^{t}p^{1/2}, \quad 0 \le t \le 1.
\]
It is the unique geodesic in $P$ from $p$ (at $t=0$) to $q$ (at $t=1$). 
It is known that $p \sharp_{t}q = q \sharp_{1-t}p$ and $(p \sharp_t q)^{-1} = p^{-1} \sharp_t q^{-1}$.
When $t = 1/2$, we abbreviate $p \sharp_{1/2} q$ as $p \sharp q$. 

Similarly, the $t$-spectral mean of $p, q \in P$ is then defined as 
\[
p \natural_t q = (p^{-1} \sharp q)^{t} p (p^{-1} \sharp q)^{t}, \quad 0 \le t \le 1.
\]
When $t = 1/2$, we abbreviate {$p \natural_{1/2} q$ as $p \natural q$.}

The following result is a generalization of \eqref{E:KL07-gm} to semisimple Lie groups.

\begin{theorem} \label{T:orbit-sum-geometric}
For $X, Y \in \p$, there exist $u, v \in K$ such that 
\begin{equation} \label{E:orbit-sum-geometric}
e^{2X}\sharp e^{2Y} = e^{\Ad u(X) + \Ad v(Y)}.
\end{equation}
\end{theorem}

\begin{proof}
Suppose $e^{2X}\sharp e^{2Y} = e^Z$ for some $Z \in \p$. 
Then by definition 
\[
e^Ze^{-2X}e^{Z} = [e^X(e^{-X}e^{2Y}e^{-X})^{1/2}e^X]e^{-2X}[e^X(e^{-X}e^{2Y}e^{-X})^{1/2}e^X] = e^{2Y}.
\] 
Now by Theorem \ref{T:orbit-sum-exponential}, there exist $k_1, k_2 \in K$ such that  
\[
e^{2Y} = e^Ze^{-2X}e^{Z} = e^{\Ad k_1(2Z) + \Ad k_2(-2X)}.
\]
Since the restriction of exponential map on $\p$ is one-to-one, we have 
\[
2Y = \Ad k_1(2Z) + \Ad k_2(-2X) =  2 \Ad k_1(Z) - 2\Ad k_2(X).
\]
That is, $\Ad k_1(Z) = \Ad k_2(X) + Y$. Hence $Z = \Ad u(X) + \Ad v(Y)$ with $u = k_1^{-1}k_2$ and $v = k_1^{-1}$.
\end{proof}

The following result is a generalization of \eqref{E:KL07-sm1} to semisimple Lie groups.

\begin{lemma} \label{L:spectral-mean}
For $p, q \in P$,  there exists a unique $k \in K$ such that 
\[
p \natural q = k(p^{1/2}qp^{1/2})^{1/2}k^{-1}.
\]
Consequently, $(p \natural q)^2$ is $K$-conjugate to $p^{1/2}qp^{1/2} \in P$ and $G$-conjugate to $pq \in G$.
\end{lemma}

\begin{proof} Let $g = (p^{-1}\sharp q)^{1/2}p^{1/2}$ and let $g = kr$ be the right Cartan decomposition of $g$ with $k \in K$ and $r \in P$. 
Then $g^*g = r^2$ and $gg^* = kr^2k^{-1} = k(g^*g)k^{-1}$. Therefore,
\begin{align*}
p \natural q &= [(p^{-1}\sharp q)^{1/2}p^{1/2}][(p^{-1}\sharp q)^{1/2}p^{1/2}]^* \\
               &= k([(p^{-1}\sharp q)^{1/2}p^{1/2}]^*[(p^{-1}\sharp q)^{1/2}p^{1/2}])k^{-1}\\
               &= k(p^{1/2}(p^{-1}\sharp q)p^{1/2})k^{-1} \\
               &= k(p^{1/2}qp^{1/2})^{1/2}k^{-1}. \qedhere
\end{align*}
\end{proof}

The following result is a generalization of \eqref{E:KL07-sm} to semisimple Lie groups.

\begin{theorem} \label{T:orbit-sum-spectral}
For $X, Y \in \p$, there exist $u, v \in K$ such that 
\begin{equation} \label{E:orbit-sum-spectral}
e^{2X}\natural e^{2Y} = e^{\Ad u(X) + \Ad v(Y)}.
\end{equation}
\end{theorem}
\begin{proof}
By Theorem \ref{T:orbit-sum-exponential}, there exist $k_1, k_2 \in K$ such that  
\[
(e^{X}e^{2Y}e^{X})^{1/2}= \left(e^{\Ad k_1(2X) + \Ad k_2(2Y)}\right)^{1/2} = e^{\Ad k_1(X) + \Ad k_2(Y)}.
\]
According to Lemma \ref{L:spectral-mean}, there exists a unique $k_3 \in K$ such that 
\[
e^{2X}\natural e^{2Y} = k_3(e^{X}e^{2Y}e^{X})^{1/2}k_3^{-1}.
\]
Combining the above equalities, we have
\begin{align*}
e^{2X}\natural e^{2Y} &= k_3(e^{X}e^{2Y}e^{X})^{1/2}k_3^{-1} \\
                               &= k_3e^{\Ad k_1(X) + \Ad k_2(Y)}k_3^{-1} \\
                               &= e^{\Ad k_3(\Ad k_1(X) + \Ad k_2(Y))} \\
                               &= e^{\Ad u(X) + \Ad v(Y)}
\end{align*}
with $u = k_3k_1$ and $v = k_3k_2$.
\end{proof}

\section{Two means for symmetric spaces}
Let the notation be as in Section 3.
An element $X \in \g$ is called real semisimple (resp., nilpotent) if $\ad X$ is diagonalizable over $\R$ (resp., nilpotent). An element $g \in G$ is called {\em hyperbolic} (resp., {\em unipotent}) if $g = \exp X$ for some real semisimple (resp., nilpotent) $X \in \g$; in either case $X$ is unique and we write $X = \log g$. An element $g \in G$ is called {\em elliptic} if $\Ad g$ is diagonalizable over $\C$ with eigenvalues of modulus $1$. According to \cite[Proposition 2.1]{Ko73}, each $g \in G$ can be uniquely written as
\begin{equation} \label{E:CMJD}
g = ehu,
\end{equation}
where $e$ is elliptic, $h$ is hyperbolic,  $u$ is unipotent, and the three elements $e, h$ and $u$ commute.
The decomposition \eqref{E:CMJD} is called the {\em complete multiplicative Jordan decomposition}, abbreviated as CMJD.

The Weyl group $W$ of $(\g, \a)$ acts simply transitively on $\a$ (and also on $A$ through the exponential map $\exp: \a \to A$).
For any real semisimple $X \in \g$, let $W(X)$ denote the set of elements in $\a$ that are conjugate to $X$, i.e.,
\[
W(X) = \Ad G(X) \cap \a.
\]
It is known from \cite[Proposition 2.4]{Ko73} that  $W(X)$ is a single $W$-orbit in $\a$.
Let conv\,$W(X)$ be the convex hull in $\a$ generated by $W(X)$. For each $g \in G$, define
\[
A(g) =\exp\mbox{conv}\, W(\log h(g)),
\]
where $h(g)$ is the hyperbolic component of $g$ in its CMJD.

Kostant's {pre-order $\prec_G$} on $G$  is defined (see \cite[p.426]{Ko73}) by {setting $f \prec_G g$} if 
\[
A(f) \subset A(g).
\]
This pre-order induces a partial order on the conjugacy classes of $G$. It is known from \cite[Theorem 3.1] {Ko73} that this {pre-order $\prec_G$} does not  depend on the choice of $\a$.

The following result  is a generalization of Theorem \ref{T:GT-ALT-Spectral} to semisimple Lie groups.

\begin{theorem}  \label{T:GT-ALT-Spectra-Gl}
Let $X, Y \in \p$. For all $r > 0$, let 
\[
\phi (r) = (e^{rX} \sharp \, e^{rY})^{2/r} \quad \text{and} \quad \psi (r) = (e^{rX} \natural \, e^{rY})^{2/r}.
\]
Then the following statements are valid.
\begin{enumerate}
\item $\phi (r)$ is monotonically decreasing on $r$ with respect to {$\prec_G$}.
\item $\psi (r)$ is monotonically increasing on $r$ with respect to {$\prec_G$}.
\item $\displaystyle \lim_{r \to 0} \phi(r) = e^{X+Y} =  \lim_{r \to 0} \psi(r)$.
\item {$(e^{rX} \sharp \, e^{rY})^{2/r} \prec_G e^{X+Y} \prec_G (e^{rX} \natural \, e^{rY})^{2/r}, \quad \mbox{for all } \, r > 0$}.
\end{enumerate}
\end{theorem}

\begin{proof}
The statements about geometric mean are special cases of Theorem 6.6, Theorem 6.7, and Theorem 6.8 in \cite{TL18}.
Lemma \ref{L:spectral-mean} and \cite[Theorem 3.30]{TL18} together yield the statements about spectral mean.
\end{proof}

The following result is a consequence of the fact that the $t$-geometric mean joining $p$ and $q$ is a geodesic under   Riemanniann metric inherited from the symmetric space $G/K$ \cite {LLT14}.

\begin{proposition}\label{sharp-r+s} Let $p, q \in P$. Then
\begin{equation}\label{eq:sharp-r+s}
p \sharp_{r+s} q
= (p \sharp_r q) \sharp_{\frac {s}{1-r}} q, \qquad \mbox{for all } \, 0 \le r, s, r+s \le 1.
\end{equation}
\end{proposition}

\begin{proof}
The path $\gamma(t) = p \sharp_t q$ with $0 \le t \le 1$ is the geodesic joining $p$ and $q$ in $P$ \cite {LLT14}, and it meets $p \sharp_{r+s}q$ when $t = r+s$. Given $r \in [0, 1]$, the path $\xi (t) = (p \sharp_r q)\sharp_t q$ with $0 \le t \le 1$ is the geodesic joining $(p \sharp_r q)$ and $q$ in $P$ and it meets $p \sharp_{r+s}q$ when $t = \frac {s}{1-r}$, since $\xi$ is part of $\gamma$.
\end{proof}

\begin{lemma}\label{lem:equ}
For $p, q \in P$, $g = p\sharp q$ if and only if $q = gp^{-1}g$. 
\end{lemma}

\begin{proof}
Since $g = p\sharp q$, by direct computation, we have 
$$
g p^{-1}g = (p\sharp q)p^{-1}(p\sharp q) 
= p^{1/2}\left(p^{-1/2}qp^{-1/2}\right)^{1/2}p^{1/2} p^{-1}p^{1/2}\left(p^{-1/2}qp^{-1/2}\right)^{1/2}p^{1/2} 
= q.
$$
If $q = gp^{-1}g$, then 
\[{p^{-1/2}qp^{-1/2} = p^{-1/2}gp^{-1}gp^{-1/2} = \left(p^{-1/2}gp^{-1/2}\right)^2.}\qedhere\]
\end{proof}

The following result extends some properties of $t$-spectral mean in Proposition \ref{T:t-natural}.
\begin{proposition} \label{T:t-natural-G}
Let $p, q \in P$ and $t \in [0, 1]$. Then
\begin{enumerate}
\item $(p \natural_t q)^{-1} = p^{-1} \natural_t q^{-1}$ and $p \natural_t q = q \natural_{1-t} p$.
\item $p^{-1} \sharp (p \natural_t q) = (q \natural_t p)^{-1} \sharp q = (p^{-1} \sharp q)^t$. \label{t-natural:2}
\item If $c_t = p^{-1} \sharp (p \natural_t q)$, then $p \natural_t q = c_t p c_t$ and $q \natural_t p ={ c_t^{-1}q c_t^{-1}}$. 
\item $(p \natural_r q) \natural_t (p \natural_s q) = p \natural_{(1-t)r+ts} q$ for all ${t, r, s \in [0,1]}$.
\end{enumerate}
\end{proposition}

\begin{proof}
\begin{enumerate}
\item Since $(p \sharp q)^{-1} = p^{-1} \sharp q^{-1}$, we have
\[
(p \natural_t q)^{-1} = (p^{-1} \sharp q)^{-t}p^{-1} (p^{-1} \sharp q)^{-t} = (p \sharp q^{-1})^{t}p^{-1} (p \sharp q^{-1})^{t} = p^{-1} \natural_t q^{-1}.
\]
Since $(q^{-1} \sharp p) q (q^{-1} \sharp p) = p$, we have
\begin{eqnarray*}
q \natural_{1-t} p &=& (q^{-1} \sharp p)^{1-t} q (q^{-1} \sharp p)^{1-t} = (q^{-1} \sharp p)^{-t}[(q^{-1} \sharp p) q (q^{-1} \sharp p)](q^{-1} \sharp p)^{-t} \\
&= &
(p^{-1} \sharp q)^t p (p^{-1} \sharp q)^t = p \natural_t q.
\end{eqnarray*}
\medskip
\item By direct computation, we have
\begin{align*}
p^{-1} \sharp (p \natural_t q) &= p^{-1/2}\left(p^{1/2}(p \natural_t q)p^{1/2}\right)^{1/2}p^{-1/2}\\
&= p^{-1/2}\left(p^{1/2}(p^{-1} \sharp q)^{t} p (p^{-1} \sharp q)^{t}p^{1/2}\right)^{1/2}p^{-1/2}\\
&= p^{-1/2}\left([p^{1/2}(p^{-1} \sharp q)^{t} p^{1/2}][p^{1/2} (p^{-1} \sharp q)^{t}p^{1/2}\right])^{1/2}p^{-1/2}\\
&= p^{-1/2}\left(p^{1/2}(p^{-1} \sharp q)^{t} p^{1/2}\right)p^{-1/2}\\
&= (p^{-1} \sharp q)^{t}.
\end{align*}
Then $(q \natural_t p)^{-1} \sharp q= q\sharp (q \natural_t p)^{-1} = (q^{-1} \sharp (q \natural_t p))^{-1} = [(q^{-1} \sharp p)^t]^{-1} = (p^{-1} \sharp q)^t$.
\medskip
\item By (2), we have $c_t = (p^{-1} \sharp q)^t$ and thus $c_t p c_t = p \natural_t q$ by definition.
Since $c_t^{-1} = (q^{-1} \sharp p)^t$, we have ${c_t^{-1}qc_t^{-1}} = q \natural_t p$.
\medskip
\item 
Let $c_i = (p^{-1} \sharp q)^i$ for all $i \in [0, 1]$. Then $p\natural_i q = c_i pc_i$ by (3). It follows that 
\begin{equation}\label{c_{i+j}}
c_{i+j} = c_ic_j,
\quad c_{ij} = c_i^j\qquad \mbox{for all } j \in [0, 1],\quad 0 \le i+j \le1.
\end{equation}
We first assume $r \le s$. Since $c_rpc_r, c_spc_s \in P$ and 
\[
c_spc_s = c_{s-r}c_rpc_rc_{s-r}= c_{s-r}[(c_rpc_r)^{-1} ]^{-1}c_{s-r},
\]
we have $(c_rpc_r)^{-1} \sharp (c_spc_s)  = c_{s-r}$ by Lemma \ref{lem:equ}.
Then
\begin{align*}
(p \natural_r q) \natural_t (p \natural_s q) &= (c_rpc_r) \natural_t (c_spc_s)\\
&= [(c_rpc_r)^{-1} \sharp (c_spc_s)]^t (c_rpc_r)[(c_rpc_r)^{-1} \sharp (c_spc_s)]^t \\
&= c_{s-r}^t (c_rpc_r)c_{s-r}^t \\
&= c_{t(s-r)} (c_rpc_r)c_{t(s-r)}  \\
&= c_{t(s-r)+r}pc_{t(s-r)+r}  \\
&= c_{(1-t)r+ts} p c_{(1-t)r+ts} \\
&= p \natural_{(1-t)r+ts} q.
\end{align*}
If $r \ge s$, then by (1) we have 
\[
(p \natural_r q) \natural_t (p \natural_s q) = (p \natural_s q) \natural_{1-t} (p \natural_r q) = p \natural_{[1-(1-t)]s + (1-t)r} q = p \natural_{(1-t)r+ts} q.\qedhere
\]
\end{enumerate}
\end{proof}

\begin{proposition}\label{natural-r+s}
Let $p, q \in P$. Then
\begin{equation}\label{eq:natural-r+s}
p \natural_{r+s} q = (p \natural_r q) \natural_{\frac {s}{1-r}} q, \qquad \mbox{for all }\, 0 \le r, s, r+s \le 1.
\end{equation}
\end{proposition}

\begin{proof}
According to Proposition \ref{T:t-natural-G} (4), we have (by assigning $s = 1$ and $t = \frac{s}{1-r}$, where the first $s$ is the $s$ in Proposition \ref{T:t-natural-G} (4) and the second $s$ is the $s$ in this proposition), 
\[
 (p \natural_r q) \natural_{\frac {s}{1-r}} q =  (p \natural_r q) \natural_t (p \natural_1 q) = p \natural_{(1-t)r+t} q = p \natural_{r+(1-r)t} q = p \natural_{r+s} q. \qedhere
\]
\end{proof}

The following figures illustrate both sides of Proposition~\ref{T:t-natural-G} (4).

\begin{minipage}{\linewidth}
      \centering
      \begin{minipage}{0.45\linewidth}
          \begin{figure}[H]
              \includegraphics[width=\linewidth]{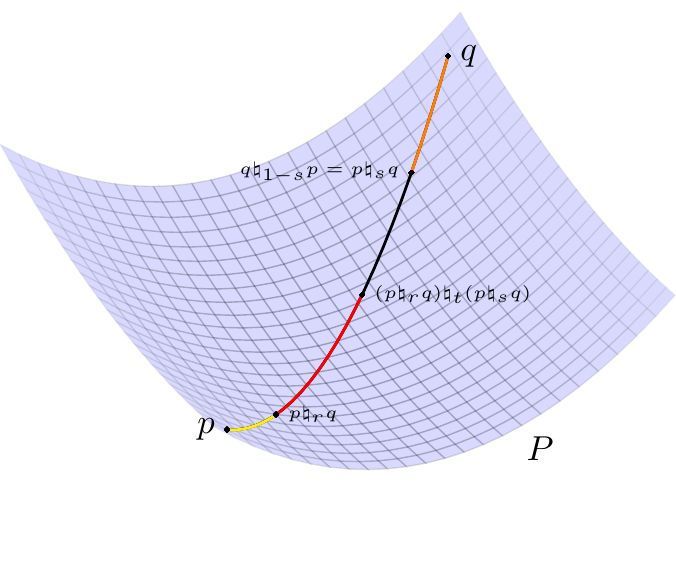}
              \caption{$(p \natural_r q) \natural_t (p \natural_s q)$}
          \end{figure}
      \end{minipage}
      \hspace{0.05\linewidth}
      \begin{minipage}{0.45\linewidth}
          \begin{figure}[H]
              \includegraphics[width=\linewidth]{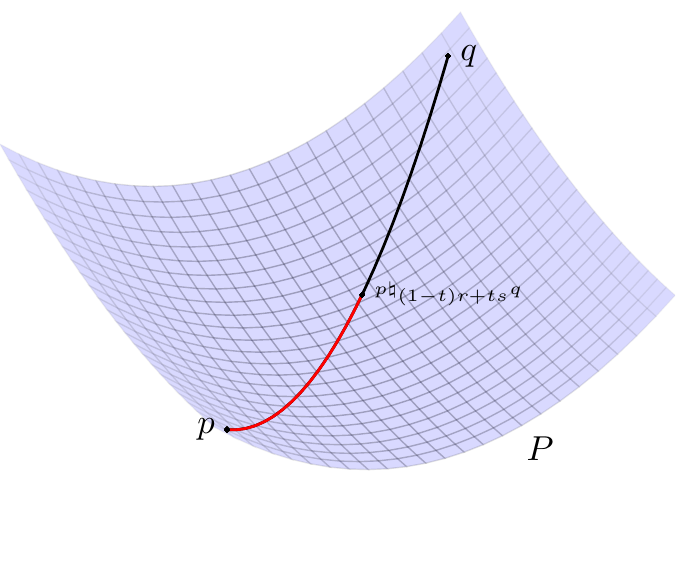}
              \caption{$p \natural_{(1-t)r+ts} q$}
          \end{figure}
      \end{minipage}
  \end{minipage}

\vspace{.5cm}

Proposition \ref{sharp-r+s} and Proposition \ref{natural-r+s} provide  information about how the curves behave. Proposition \ref{sharp-r+s} is a consequence of the fact that the $t$-geometric mean joining $p$ and $q$ is a geodesic.

Finally, we extend Theorem \ref{T:log-majorization} to symmetric space by using Theorem \ref{T:log-majorization} itself and applying Kostant's theorem as in the proof of  \cite [Theorem 3.5]{LLT14}.

\begin{theorem} \label{T:sharp-natural-G}
Let $p, q \in P$. Then
\[
	{p \sharp_t q \prec_G p \natural_t q}, \quad  \mbox{for all } \, 0 \le t \le 1.
\]
\end{theorem}

\begin{proof}

A result of Kostant \cite[Theorem 3.1]{Ko73} asserts that for any given
$f, g \in G$, {$f\prec_G g$} if and only if $|\pi(f)|\le |\pi(g)|$ for any finite dimensional representation $\pi$ of $G$, where $|\pi(g)|$ is the spectral radius of $\pi(g)$.
Let $\pi: G\to \GL(V)$ be any finite dimensional representation of $G$.
There exists an inner product on $V$ such that $\pi(z)$ is positive definite for all $z\in P$ (see \cite[p.435]{Ko73}).  
By Kostant's result, it suffices to show that $|\pi(p \sharp_t q)| \le |\pi(p\natural_t q)|$.
Then we have
\begin{eqnarray*}
	|\pi(p \sharp_t q)|&=& |\pi(p^{1/2}(p^{-1/2}q p^{-1/2})^t p^{1/2})| \\
	&=& |\pi(p)^{1/2}(\pi(p)^{-1/2}\pi(q)\pi(p)^{-1/2})^t\pi(p)^{1/2}| \\
	&=& |\pi(p)\sharp_t \pi(q)|\\
	&\le & |\pi(p)\natural_t \pi(q)| \qquad \text{(by Theorem~\ref{T:log-majorization})}\\
	&=& |(\pi(p)^{-1}\sharp \pi(q))^t \pi(p) (\pi(p)^{-1}\sharp \pi(q))^t|\\
	&=& |\pi((p^{-1}\sharp q)^t p (p^{-1}\sharp q)^t)|\\
	&=& |\pi(p\natural_t q)|.
\end{eqnarray*}
\end{proof}

\noindent{\bf Acknowledgement }
{We are thankful to the anonymous referee for the careful reading of our paper and for the
constructive suggestions that helped us to make improvement. In particular, we are grateful to the referee for bringing our attention to the two references Lee and Lim~\cite{LL07} and Kim~\cite{Ki21}.}
%%%%%%%%%%%%%%%%%%%%%%%%%%%%%%%%%%%%%%%%%%%%%%%%%%%%%%%
\bibliographystyle{abbrv}
\bibliography{refs}

\end{document}